\documentclass[12pt,a4paper]{article}
\usepackage{fancyhdr}
\usepackage[latin1]{inputenc}
\usepackage{amsmath,amssymb,amsthm,amsfonts}
\usepackage{graphicx}
\usepackage{color}
\usepackage{mathrsfs,mathptmx}
\usepackage{graphs,graphicx,xcolor}
\usepackage{fullpage}
\usepackage[T1]{fontenc}
\newcommand{\N}{{\mathbb N}}

\newcommand{\Z}{{\mathbb Z}}
\newcommand{\uno}{{\mathbf{1}}}

\newtheorem{theorem}{Theorem}
\newtheorem{lemma}[theorem]{Lemma}
\newtheorem{corollary}[theorem]{Corollary}
\newtheorem{conjecture}[theorem]{Conjecture}

\newtheorem{remark}[theorem]{Remark}

\begin{document}

\begin{verbatim}\end{verbatim}\vspace{2.5cm}

\begin{center}
\textsc{\Large\textbf{The Minor inequalities in the description of the Set Covering Polyhedron of Circulant Matrices}}
\end{center}
\smallskip

\begin{center}
\large Silvia M.\,Bianchi\footnote{sbianchi@fceia.unr.edu.ar} \;\;\; Graciela L.\,Nasini\footnote{nasini@fceia.unr.edu.ar} \;\;\; Paola B.\,Tolomei\footnote{ptolomei@fceia.unr.edu.ar}
\end{center}

\smallskip

\begin{center}
\textit{Departamento de Matem\'atica, Facultad de Ciencias Exactas, Ingenier\'ia y Agrimensura, Universidad Nacional de Rosario, 2000 Rosario, Santa Fe, Argentina}
\medskip

\textit{and CONICET, Argentina}
\end{center}

\medskip
\begin{abstract}
In this work we give a complete description of the set covering polyhedron of circulant matrices $C^k_{sk}$ with $s = 2,3$ and $k \geq 3 $  by linear inequalities. In particular, we prove that every non boolean facet defining inequality is associated with a circulant minor of the matrix.  
We also give a polynomial time separation algorithm for inequalities involved in the description. 

\medskip

\end{abstract}

\section{Introduction}
\label{sect:introduction}

The \emph{weighted set covering problem} can be stated as 
\[
\begin{array}{lrclr}
\min &  c^Tx &      &    & \\
     &   A x & \geq & \mathbf 1 & \ \ \ \ \ \ \mbox{(SC)} \\  
     &     x & \in  & \{0,1\}^n & \\
\end{array}
\]
where $A$ is an $m\times n$ matrix with $0,1$ entries, $c$ is an $n$-vector and  $\mathbf 1$ is the $m$-vector of all ones.

These types of problems are relevant in practice, but hard to solve in general. One often successful way to tackle such problems is the polyhedral approach involving the solution space of the problem \cite{GLS}. 

The \emph{set covering polyhedron} (SCP), denoted by $Q(A)$, is the convex hull of integer 
points in $Q_R(A)=\{x\in [0,1]^n: Ax\geq \mathbf 1 \}$.

If for some matrix $A$ it holds that $Q_R(A)=Q(A)$, the matrix is called \emph{ideal} and this would enable us to solve SC as a linear program using the constraints $x\geq \mathbf 0$ instead of the integrality requirements.
However, when $A$ is nonideal, finding  a description for $Q(A)$ in terms of linear restrictions is, in general, as hard as solving SC.

In \cite{CN} Cornu\'ejols and Novick studied the SCP 
on  a particular class of matrices, called \textit{circulant matrices} and denoted as $C_n^k$ with $1\leq k\leq n-1$. 
They identified all the ideal circulant matrices which are 
$C_6^3,C_9^3, C_8^4$ and $C_n^2$,  for even  $n\geq 4$. 
They also provide sufficient conditions a given submatrix must satisfy in order to be a \emph{circulant minor}. Circulant matrices and circulant minors will be formally defined in the next section.

Using these results, Argiroffo and Bianchi obtained in \cite{circu} a family of facets of $Q(C_n^k)$ associated with some particular circulant minors.
Previously, in \cite{Mah} Bouchakour et al., when working on the dominating set polyhedron of cycles, they obtained the complete description of $Q(C_n^3)$ for every $n\geq 5$. Interestingly, all the non boolean non rank constraints involved in this description belong to this family of inequalities associated with circulant minors.

Later, Aguilera in \cite{Nes} completely identified all circulant minors that a circulant matrix may have. This result allowed us to obtain in \cite{BNT} a wider class of valid inequalities associated with circulant minors, which we  
call \emph{minor inequalities}. 

In this paper, we present two new families of circulant matrices for which the SCP can be described by boolean facets and minor inequalities. We also give a polynomial time separation algorithm for these inequalities. 

A preliminary version of this work appeared without proofs in \cite{tunez}.   

\section{Notations, definitions and preliminary results}

Given a $0,1$ matrix $A$, we say that a row $v$ of $A$ is a \textit{dominating row} if $v\geq u$ for some $u$ row of $A$, $u\neq v$. 
In this work, every matrix has $0,1$ entries, no zero columns and no dominating rows. 

Also, every time we state $S\subseteq \Z_n$ for some $n\in \N$, we consider $S\subseteq \{0,\dots,n-1\}$ and the addition between the elements of $S$ is taken modulo $n$.  
Rows and columns of an $m\times n$ matrix $A$ are indexed by $\Z_m$ and $\Z_n$ respectively.
Two matrices $A$ and $A'$ are \emph{isomorphic} 
if $A'$ can be obtained from $A$ by permutation of rows and columns.

Given $N\subset \Z_n$, the \emph{minor of} $A$ \emph{obtained by contraction of} $N$ and denoted by $A/N$, is the submatrix of $A$ that results after removing all columns indexed in $N$ and all the dominating rows that may occur.
In this work, when we refer to \emph{a minor} of $A$ we always consider a minor obtained by contraction. 

Considering the one-to-one correspondence between a vector $x\in \{0,1\}^n$ and the subset $S_x\subseteq \Z_n$ whose characteristic vector is $x$ itself, we agree to abuse of notation by writing $x$ instead of $S_x$.

A \emph{cover} of a matrix $A$ is a vector $x\in \{0,1\}^n$ such that $Ax\geq \uno$. A cover $x$ of $A$ is \textit{minimal} if there is no other cover $\tilde{x}$ such that $\tilde{x}\subset x$.
A cover $x$ is \emph{minimum} if $|x|=\sum_{i\in\Z_n} x_i$ is minimum and in this case $|x|$ is called the \emph{covering number} of the matrix $A$, denoted by $\tau(A)$.
Since every cover of a minor of $A$ is a cover of $A$, it holds that $\tau(A/N)\geq \tau(A)$ for all $N\subset \Z_n$.

\medskip

We denote with $(a,b)_n$ the
$\Z_n$-cyclic open interval of points strictly between $a$ and $b$ and 
analogous meanings for $[a,b)_n$, $(a,b]_n$ and $[a,b]_n$.

Given $n$ and $k$ with $1\leq k\leq n-1$, the \emph{circulant} matrix $C_n^{k}$ is the square matrix whose $i$-th row is the incidence vector of $C^{i}=[i,i+k)_n$.

It is not hard to see that, for every $i \in \Z_n$,  
$$x^i=\left\{i+hk : 0\leq h \leq \left\lceil \frac{n}{k}\right\rceil-1 \right\}\subset \Z_n$$ 
is a cover of $C^k_n$ of size $\left\lceil \frac{n}{k}\right\rceil$. It is also clear that $\tau(C_n^k)\geq\left\lceil \frac{n}{k}\right\rceil$ and then $\tau(C_n^k)=\left\lceil \frac{n}{k}\right\rceil$.
Let us also observe that for every minimal cover $x$ of $C^k_n$ and any $i\in \Z_n$, $|x\cap C^i|\leq 2$.

We say that a minor of $C_n^k$ is a \textit{circulant minor} if it is isomorphic to a circulant matrix.
In \cite{CN}, the authors give sufficient conditions for a subset $N\subset \Z_n$ to ensure that $C_n^k/N$ is a circulant minor of $C^k_n$. These conditions are obtained in terms of simple dicycles in a particular digraph. 

Indeed, given $C^{k}_n$, the digraph $G(C^{k}_n)$ has vertex set $\Z_n$ and $(i,j)$ is an arc of $G(C^{k}_n)$ if $j\in \{i+k,i+k+1\}$. We say that an arc $(i,i+k)$ has length $k$ and an arc $(i,i+k+1)$ has length $k+1$. 

If $D$ is a simple dicycle in $G(C^{k}_n)$, and $n_2$ and $n_3$ denote the number of arcs of length $k$ and $k+1$ respectively, $k n_2 +(k+1)n_3=n_1 n$ for some unique positive integer $n_1$.
We say that $n_1, n_2$ and $n_3$ are the \emph{parameters associated with} the dicycle $D$.

In \cite{Nes2} it is proved that the existence of nonnegative integers $n_1, n_2$ and $n_3$ satisfying the conditions  $n_1 n = k n_2 + (k+1) n_3$ and $\gcd(n_1,n_2,n_3)=1$ are also sufficient for the existence of a simple dicycle in $G(C^{k}_n)$ with $n_2$ arcs of length $k$ and $n_3$ arcs of length $k+1$. Moreover, the same author completely characterized in \cite{Nes} subsets $N$ of $\Z_n$ for which $C_n^k/N$ is a circulant minor  in terms of dicycles in the digraph $G(C^{k}_n)$. We rewrite Theorem 3.10 of \cite{Nes} in order to suit the current notation in the following way:

\begin{theorem} \label{minorsgrafo}
Let $n,k$ be integers verifying $2\leq k \leq n-1$ and let $N\subset \Z_n$ such that $1\leq |N|\leq n-2$.
Then, the following are equivalent:
\begin{enumerate}
\item $C_n^k/N$ is isomorphic to $C_{n'}^{k'}$.  
\item $N$ induces in $G(C^{k}_n)$ $d\geq 1$ disjoint simple dicycles $D_0,\ldots,D_{d-1}$, each of them having the same parameters $n_1$, $n_2$ and $n_3$ and such that 
$|N|=d(n_2+n_3)$, 
$n'=n-d(n_2+n_3)\geq 1$ and $k'=k-dn_1\geq 1$.
\end{enumerate}
\end{theorem}

Thus, whenever we refer to a circulant minor of $C_n^k$ with parameters $d$, $n_1$, $n_2$ and $n_3$, we mean the non negative integers whose existence is guaranteed by the previous theorem. 
In addition, for each $j\in \Z_d$, $N^j$ refers to the subset of $\Z_n$ inducing the simple dicycle $D^j$ in $G(C^k_n)$, $W^j=\{i\in N^j: i-(k+1)\in N^j\}$ and $W=\cup_{j\in\Z_d} W^j$.   
Then, for all $j\in \Z_d$, $\left|W^j\right|=n_3$ and $\left|N^j\right|=n_2+n_3$.

In \cite{BNTejor} it was proved that 
given $W\subset \Z_n$ corresponding to a circulant minor of $C^k_n$, we can rebuilt $N$ such that $C_n^k/N$ is such a minor.
So, in what follows, we usually refer to a circulant minor \emph{defined by} $W$. We also say that $W$ \emph{defines} a circulant minor.

Circulant minors, or equivalently subsets $W\subset \Z_n$ inducing them, play an important role in
the description of the set covering polytope of circulant matrices.

It is known that  $Q(C_n^k)$ is a full dimensional polyhedron. Also, for every $i\in \Z_n$, the constraints $x_i\geq 0$, $x_i\leq 1$ and $\sum_{j\in C^{i}} x_j \geq 1$ are facet defining inequalities of $Q(C_n^k)$ and we call them \emph{boolean facets}.  
In addition, it is known that every non boolean facet of $Q(C_n^k)$ has positive coefficients \cite{circu}.

The inequality  $\sum_{i=1}^n x_i\geq \left\lceil \frac{n}{k}\right\rceil$, called the \emph{rank constraint}, is always valid for $Q(C_n^k)$ and defines a facet if and only if $n$ is not a multiple of $k$ (see \cite{Sa}).

However, for most circulant matrices these constraints are not enough to obtain their corresponding SCP \cite{circu}.
Actually, in \cite{BNT} the authors obtained a new family of non boolean no rank facets defining inequalities of the SCP of circulant matrices associated with circulant minors.

\begin{lemma} \cite{BNT} \label{nosotras}
Let $W\subset\Z_n$ define a circulant minor of $C_n^k$ isomorphic to $C_{n'}^{k'}$. Then, the  inequality 
\begin{equation}\label{ecu}
\sum_{i\in W} 2 x_i + \sum_{i\notin W} x_i \geq \left\lceil
\frac{n'}{k'}\right\rceil
\end{equation}
is a valid inequality for $Q(C_n^k)$.
Moreover, if $2\leq k'\leq n'-2$, $\left\lceil
\frac{n'}{k'}\right\rceil > \left\lceil
\frac{n}{k}\right\rceil$  and $n'= 1(\mathrm{mod} \,k')$ then the inequality (\ref{ecu}) defines a facet of $Q(C^k_n)$.
\end{lemma}

From now on, we say that inequality (\ref{ecu}) is the \emph{minor inequality corresponding to} $W$ or to the minor $W$ defines.

Observe that if $\left\lceil\frac{n'}{k'}\right\rceil=\left\lceil\frac{n}{k}\right\rceil$, the minor inequality is dominated by the rank constraint.
In addition, in \cite{BNT} it is proved that if $n'$ is a multiple of $k'$ then the corresponding inequality is 
valid for $Q_R(C_n^k)$. 
As our main interest are the relevant constraints in the description of $Q(C_n^k)$, we call \emph{relevant minors} to those minors isomorphic to $C_{n'}^{k'}$ with $n'\neq 0\, (\mathrm{mod}\, k')$ and $\left\lceil\frac{n'}{k'}\right\rceil>\left\lceil\frac{n}{k}\right\rceil$. Inequalities associated with relevant minors will be \emph{relevant minor inequalities}.

In \cite{BNT} we stated the following conjecture:

\begin{conjecture} \cite{BNT} \label{conj2}
A relevant minor inequality corresponding to a minor of $C_n^k$ isomorphic to $C^{k'}_{n'}$ defines a facet of $Q(C^k_n)$ if and only if $n'=1\, (\mathrm{mod} \,k')$. 
\end{conjecture}

It can be seen that every non boolean facet defining inequality of $Q(C^3_n)$ obtained in \cite{Mah} is either the rank constraint or it is a relevant minor inequality satisfying Conjecture \ref{conj2}.
Our goal is to enlarge the family of circulant matrices for which the same holds.
For this pourpose in section 3 we obtain necessary conditions for 
an inequality to be a 
non boolean non rank facet defining inequality of $Q(C^k_{n})$. In section 4, we focus on matrices of the form $C_{sk}^k$ and find that every facet defining inequality with right hand side $s+1$ is a minor inequality satisfying Conjecture \ref{conj2}. Moreover, we prove that this inequalities can be separated in polytime. Finally, in section 5 we prove that $Q(C_{2k}^k)$ and $Q(C_{3k}^k)$ are described in terms of boolean facets and minor inequalities with right hand side $s+1$.

\section{Properties of facets of $Q(C^k_{n})$} \label{propiedadesfacetas}
 
Let $ax\geq \alpha$ be a non boolean, non rank facet defining inequality  of $Q(C_n^k)$ with integer coefficients. 
A root $\tilde{x}$ of $ax\geq \alpha$ is a cover of $C_n^k$ satisfying $a\tilde{x}= \alpha$.
Since the inequality has positive coefficients then $\tilde{x}$ is a minimal cover of $C_n^k$. 

We define $a^0=\min \{a_i: i\in \Z_n\}$ and $W=\{i\in \Z_n: a_i>a^0\}$.
Clearly, $a^0\geq 1$ and $W\neq \Z_n$. Moreover, $W\neq \emptyset$ since otherwise $ax\geq \alpha$ would be dominated by the rank inequality. By denoting $\overline{W}=\{i\in \Z_n: i\notin W\}$, $ax\geq \alpha$   
can be written as 

\begin{equation}\label{desig}
\sum_{i\in W} a_i x_i+a^0\sum_{i\in \overline{W}} x_i\geq \alpha. 
\end{equation}

Observe that, for every cover $\tilde{x}$ of $C^k_n$, 

$$a \tilde{x}= \sum_{i\in \tilde{x} \cap W} a_i +a^0 | \tilde{x}\cap  \overline{W}| \geq a^0|\tilde{x}|.$$

Since (\ref{desig}) is not the rank inequality, it has a root $\tilde{x}$ that is not a minimum cover and then  
$$\alpha=a \tilde{x} \geq a^0|\tilde{x}|\geq a^0(\tau(C_n^k) +1)= a^0\left(\left\lceil \frac{n}{k}\right\rceil+1\right), $$
and every minimum cover $x$ must satisfy $x\cap W \neq \emptyset$ since otherwise, it would violate (\ref{desig}). 

For the sequel it is convenient to make the next observation:

\begin{remark}\label{useful}
Every non boolean non rank facet defining inequality of $Q(C^k_n)$ is of the form (\ref{desig}) with $a^0\geq 1$, $\emptyset \subsetneq W\subsetneq \Z_n$, $a_i\geq a^0+1$ for all $i\in W$ and $\alpha \geq a^0(\left\lceil \frac{n}{k}\right\rceil+1)$. Moreover, $|x\cap W |\geq 1$ for every minimum cover $x$ of $C_n^k$.
\end{remark}

We have the following results: 

\begin{lemma}\label{prop} 
Let (\ref{desig}) be a non boolean non rank facet defining inequality of $Q(C_n^k)$.  
\begin{enumerate}
\item 
For every $i\in \Z_n$ there exists \label{item3}
\begin{enumerate}
\item
a root $\tilde{x}$ such that $i\in \tilde{x}$, 
\item
a root $\tilde{x}$ such that $i \notin \tilde{x}$, 
\item \label{tirita2}
a root $\tilde{x}$ such that $|\tilde{x}\cap C^{i}|= 2$.\label{dos}
\end{enumerate}
\item
Let $i\in W$ and $\tilde{x}$ a root such that $i\in \tilde{x}$. 
\begin{enumerate}
\item
If there exists $j\neq i$ such that $j\in \tilde{x}\cap C^{i-k+1}$, then $[i,j+k]_n \subset W$.  
\item
If there exists $j\neq i$ such that $j\in \tilde{x}\cap C^i$, then $[j-k, i]_n \subset W$.  
\end{enumerate}
\end{enumerate}
 
\end{lemma}

\proof
Let $i\in \Z_n$.

If for every root $\tilde{x}$ of (\ref{desig}) it holds that $i\notin \tilde{x}$ ($i\in \tilde{x}$) then every root of (\ref{desig}) is also a root of the boolean facet defined by $x_i\geq 0$ ($x_i\leq 1$), a contradiction. Then, items 1.(a) and (b) hold. 

Let us observe that every root $\tilde{x}$ of (\ref{desig}) satisfying $|\tilde{x}\cap C^{i}|=1$ is also a root of the boolean facet defined by the inequality $\sum_{j\in C^{i}} x_j\geq 1$. 
Then, we conclude that there exists a root $\tilde{x}$ such that $|\tilde{x}\cap C^{i}|\geq 2$. Recalling that $\tilde{x}$ is a minimal cover, $|\tilde{x}\cap C^{i}|\leq 2$ and then item 1.(c) follows.

\smallskip

In order to prove item 2., let $i\in W$ (i.e. $a_{i}>a^0$) and $\tilde{x}$ be a root of (\ref{desig}) such that $i\in \tilde{x}$.

Assume $j\neq i$ such that $j\in \tilde{x}\cap C^{i-k+1}$. Observe that for any $h\in [i,j+k]_n$, $\hat{x}=\tilde{x}\setminus \{i\}\cup \{h\}$ is a cover of $C^k_n$ satisfying (\ref{desig}). Then, we have
$$a\hat{x}=a\tilde{x}-a_{i}+a_{h}=\alpha-a_{i}+a_{h} \geq \alpha$$ 
implying $a_h\geq a_{i}>a^0$, i.e. $h\in W$.

Now, using similar arguments when $j\neq i$ such that $j\in \tilde{x}\cap C^i$ we arrive to $a_{h}\geq a_i>a^0$ for $h\in [j-k, i]_n$ and the lemma follows. 
\qed

\medskip
\medskip

From the previous results, we obtain the following relevant properties of facet defining inequalities of $Q(C^k_n)$.

\medskip

\begin{theorem}\label{wraya}
Let (\ref{desig}) be a non boolean non rank facet defining inequality of $Q(C_n^k)$. Then, 

\begin{enumerate}
\item
for every $i\in \Z_n$, $|C^{i}\cap \overline{W}|\geq 2$.  
\item
for every $i\in W$, $a_i\leq 2a^0$.
\end{enumerate}
\end{theorem}

\proof

Suppose that $|C^{i}\cap \overline{W}|=0$ for some $i\in \Z_n$. 
Let consider a root $x$ of (\ref{desig}) such that $|x\cap C^{i}|=2$ that exists according to Lemma \ref{prop} item 1.(c). From Lemma \ref{prop} item 2.(a) we get $i+k\in W$ and then $|C^{i+1}\cap \overline{W}|=0$. Iteratively using the same argument, we arrive to $W=\Z_n$, a contradiction. Thus, we have proved that  $|C^{i}\cap \overline{W}|\geq 1$ for all $i\in \Z_n$.

Now suppose that $C^{i}\cap \overline{W}=\{h\}$ for some $i\in \Z_n$.
Let $x$ be again a root of (\ref{desig}) such that $x\cap C^{i}=\{s,s'\}$ with $s\in [i,s')_n$.
If $h \in [i,s]_n$ then $s' \in W$ and by Lemma \ref{prop} item 2.(a),  $C^{s+1}\subset W$. Similarly, if $h\in [s',i+k-1]_n$, $s\in W$ and $C^{s'-k}\subset W$. In both cases, we obtain a contradiction with $|C^{i}\cap \overline{W}|\geq 1$ for all $i\in \Z_n$ as we have already noted. We conclude that $h\in(s,s')_n$. 

In particular,  we have proved that if $C^{i}\cap \overline{W}=\{h\}$ for some $i\in\Z_n$, then $i\neq h$.

Moreover, since $h\in(s,s')_n$, $s'\in W$. Applying Lemma \ref{prop} item 2.(a) we have that $[s',s+k]_n\subset W$ and then $C^{i+1} \cap \overline{W}=\{h\}$. Iteratively using the same argument,  we arrive to $C^{h}\cap \overline{W}=\{h\}$
contradicting the previous observation. 

Then, $|C^{i}\cap \overline{W}| \geq 2$ for all $i\in \Z_n$.

In order to prove item 2., let $i\in W$ and $\tilde{x}$ be a root of (\ref{desig}) such that $i\in \tilde{x}$. After item 1., $C^{i-k}\cap \overline{W}\neq \emptyset$ and there is $h\in C^{i-k}\cap \overline{W}$ such that $(h,i]_n\subset W$. Let $\ell \in C^{h}\cap \overline W$ then it holds that $i\in (h,\ell)_n$. 
Then, $\hat{x}=\tilde{x}\setminus \{i\} \cup \{h,\ell\}$ is a cover of $C^k_n$ and therefore, $a\hat{x}=a\tilde{x}-a_i+2 a^0=\alpha- a_i + 2a^0 \geq \alpha $ or equivalently, $a_i \leq 2a^0$.
\qed

\section{Minor inequalities of $Q(C_{sk}^k)$}

In this section we work with circulant matrices $C_n^k$ for which $n$ is a multiple of $k$. Firstly, we will see that almost every circulant matrix can be thought as a minor of such a matrix. 

In Theorem 2.10 of \cite{Nes} it was proved that a  matrix $C_n^k$ has a minor isomorphic to $C_{n'}^{k'}$ if and only if 

\begin{equation}
\label{cond_minor}
\frac{k'}{k}\leq \frac{n'}{n}\leq \frac{k'+1}{k+1}.
\end{equation}

As a consequence we have:

\begin{lemma}\label{menor_circ}
Let $n'=hk'+r$ with $1\leq r\leq k'-1$. Then, there exist $s$ and $k$ such that $C_{sk}^{k}$ has a minor isomorphic to $C_{n'}^{k'}$ if and only if $r\leq h-1$. 
\end{lemma}

\proof
Let $s$ and $k$ be such that $C_{sk}^{k}$ has a minor isomorphic to $C_{hk'+r}^{k'}$. Then, according to (\ref{cond_minor}), we have   
$\frac{k'}{k}\leq \frac{hk'+r}{sk }$ and it holds that $s\leq h$.

In addition, $\frac{hk'+r}{sk }\leq \frac{k'+1}{k+1}$, is equivalent to
\[n'= hk'+r \leq [s(k'+1)-(hk'+r)] k\]
Since $n'\geq 1$, it holds that $s(k'+1)-(hk'+r) >0$ and then, $s>\frac{hk'+r}{k'+1}$.

In summary,  we have that 
\[\frac{hk'+r}{k'+1}<s\leq h\]
and then $r\leq h-1$. 

Conversely, if $r\leq h-1$, it is easy to see that by taking $s=h$ and $k\geq \frac{hk'+r}{h-r}$ the condition (\ref{cond_minor}) holds and by Theorem 2.10 of \cite{Nes}, $C_{sk}^k$ has a minor isomorphic to $C_{sk'+r}^{k'}$.
\qed

\bigskip

Hence, for a fixed $k'$, 
except for a finite number of values of $n'$, matrix $C_{n'}^{k'}$ is isomorphic to a minor of some matrix $C^k_{sk}$.

Let us start the study of polyhedra $Q(C^k_{sk})$, for $s\geq 2$. 
Observe that matrices $C_{3s}^3$ have already been studied in \cite{Mah}. Moreover, for $k=4$ we take $s\geq 3$ since $Q(C_{8}^4)$ is described by boolean inequalities (see \cite{CN}). 

Remind that if $\tilde{x}$ is a minimum cover of $C^k_{sk}$, $\tilde{x}=x^j=\{j+rk, r\in \Z_s\}$ for some $j\in\Z_{sk}$.
Therefore, $x^i=x^{i+k}$ for any $i\in\Z_{sk}$ and there are exactly $k$ minimum covers ($x^i$ with $i\in \Z_k$) defining a partition of $\Z_{sk}$.

From now on, we consider facet defining inequalities of $Q(C^k_{sk})$ in the form
\begin{equation}\label{desigs+1}
\sum_{i\in W} a_ix_i+a^0\sum_{i\in \overline{W}} x_i\geq  (s+1)a^0 
\end{equation}
for some $a^0\geq 1$, $\emptyset \subsetneq W\subsetneq \Z_{sk}$ and   
$2 a^0\geq a_i\geq a^0+1$ for every $i\in W$. 
Recall that, from the results in Theorem \ref{wraya} we also know that $|\overline{W}\cap C^{i}|\geq 2$.

Inequalities in the form (\ref{desigs+1}) will be referred as $(s+1)$-\emph{inequalities} of $Q(C_{sk}^k)$.

We will prove that every facet defining $(s+1)$-inequality of $Q(C^k_{sk})$ is a minor inequality.

Observe that every root $\tilde{x}$ of (\ref{desigs+1}) has cardinality $s$ or $s+1$. 
Moreover, $\tilde{x}\cap W \neq \emptyset$ if and only if $\tilde{x}$ is a minimum cover of $C^k_{sk}$. Thus, if $i\in \tilde{x}\cap W$, $\tilde{x}=x^i$.

Hence, 

\begin{theorem} \label{2a0} Let (\ref{desigs+1}) be a facet defining $(s+1)$-inequality of $Q(C^k_{sk})$. Then, for every $i \in \Z_{sk}$, it holds:
\begin{enumerate}
\item \label{unosolo}
$|x^i\cap W|=1$. Moreover, $|W|=k$.
\item
If $i\in W$, $a_i=2 a^0$.
\end{enumerate}
\end{theorem}

\proof
Let $i\in\Z_{sk}$.
From the previous observations $|x^i \cap W |\geq 1$.

Let assume that there exist $j\neq \ell$ such that $ \{j, \ell\}\subset x^i \cap W$. Then $x^i=x^j= x^{\ell}$. So, given a root $x$ of (\ref{desigs+1}), $x$ contains $j$ if and only if it also contains $\ell$. Therefore 
every root of the inequality (\ref{desigs+1}) lies in the hyperplane $x_{\ell}- x_j=0$, a contradiction.

Then, $|x^i\cap W|=1$ and recalling that $\{x^i:i\in\Z_k\}$ defines a partition of $\Z_{sk}$, it holds that $|W|=k$.

Hence, if $i\in W$, $x^i\cap W=\{i\}$ and $x^i$ is a root of (\ref{desigs+1}). As a consequence, we have:
$$a x^i= a_i+ (s-1) a^0 \geq (s+1) a^0$$ 
and then, $a_i\geq 2 a^0$. By Theorem \ref{wraya} item 2., we have that $a_i=2a^0$.\qed

\medskip

We have proved that every facet defining $(s+1)$-inequality of $Q(C^k_{sk})$ can be written as 
\begin{equation}
2 \sum_{i\in W}  x_i + \sum_{i\in \overline{W}} x_i \geq s+1 
\label{minor s+1}
\end{equation}
where $W$ verifies $|W\cap x^i|=1$ and $|\overline{W}\cap C^{i}|\geq 2$, for every $i\in \Z_{sk}$. 

The next theorem proves that the class of minor inequalities of $Q(C^k_{sk})$ includes facet defining $(s+1)$-inequalities.

\begin{theorem}\label{minor}
Let $W\subset \Z_{sk}$ such that $|W \cap x^j|=1$ for every $j\in \Z_{sk}$. Then, the inequality (\ref{minor s+1}) is valid for $Q(C^k_{sk})$. Moreover, if $|\overline{W}\cap C^{j}|\geq 1$ for every $j\in \Z_{sk}$, it is a minor inequality.
\end{theorem}

\proof
Observe that, if $|\overline{W}\cap C^{i}|=0$ for some $i\in \Z_{sk}$, then inequality (\ref{minor s+1}) is the sum of the rank constraint and the boolean facet $\sum_{j\in C^i} x_j\geq 1$.

For the other cases, we have to prove that $W$ defines a circulant minor of $C^k_{sk}$. In particular, we will prove that $W$ corresponds to a simple dicycle $D$ in $G(C_{sk}^k)$ which verifies item 2. of Theorem \ref{minorsgrafo}. 

In order to obtain $D$ we proceed in the following way.
As $|W \cap x^j|=1$ for every $j\in \Z_k$, let $i_j\in \Z_{sk}$ such that  $x^j\cap W=\{i_j\}$ and let $t_j$ such that $0\leq t_j \leq s-1$ and  $i_j=j+t_jk$. 
For every $j\in \Z_k$, let define the $i_j i_{j+1}$-dipath $P_j$ in the digraph $G(C_{sk}^k)$ induced by the set $V_j\cup \{i_{j+1}\}$, where  
$$V_j=\{i_{j}+rk: 0\leq r\leq n_2^j\}\subset x^j, $$
and $n_2^j$ is defined according to the following cases.

If $j\neq k-1$ then 

$$
n_2^{j\;\;}=\left\{\begin{array}{ll}
(s-1)-(t_j-t_{j+1}) \,&\,\;\mathrm{ if } \;\;t_j-t_{j+1}\geq 0\\

(t_{j+1}-t_j)-1 \,&\,\;\mathrm{ otherwise }
\end{array}
\right.
$$

\noindent else
$$n_2^{k-1}=\left\{\begin{array}{lcl}
t_{0}-t_{k-1}-2 \;&\,\;\mathrm{ if } \;&\;t_{k-1}-t_0\leq -2\\

(s-2)-(t_{k-1}-t_{0}) \;&\,\;\mathrm{ if } \;&\;-1\leq t_{k-1}-t_0\leq s-2\\

s-1 \:&\,\;\mathrm{ if } \;&\;t_{k-1}-t_0= s-1.
\end{array}
\right.$$
\medskip

\begin{figure}[h]
	\centering
		\includegraphics[width=0.8\textwidth]{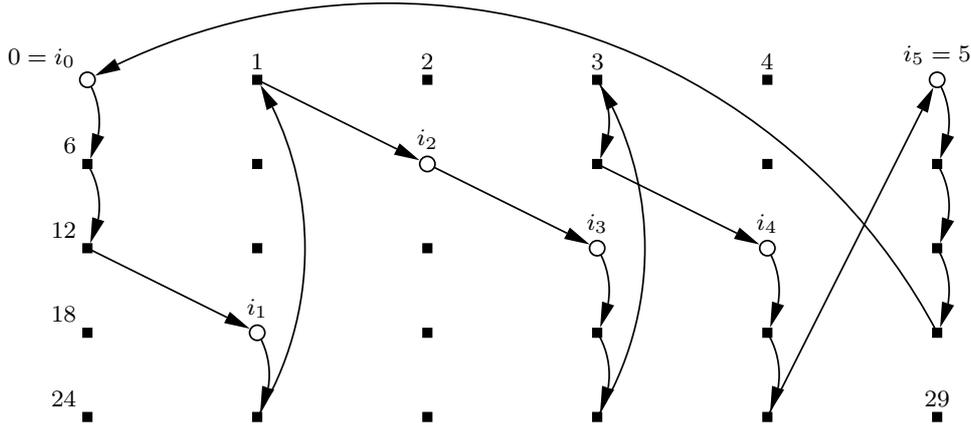}
		\caption{The dicycle $D$ associated with $W=\{0,5,8,15,16,19\}\subset \Z_{sk}$ with $s=5$ and $k=6$.}
	\label{fig:circuito}
\end{figure}

Observe that subsets $V_j$ with $j\in \Z_k$ are mutually disjoint 
and for  any $j\in \Z_k$, 
the $i_j i_{j+1}$-dipath $P_j$ induced by $V_j\cup \{i_{j+1}\}$ in the digraph $G(C_{sk}^k)$ has exactly $n_2^j$ arcs of length $k$ and one arc of length $k+1$. Moreover, $\bigcup_{j\in\Z_k} V_j$ induces a simple dicycle $D$ in $G(C^k_{sk})$ (see Figure \ref{fig:circuito} as example).

\medskip

Let us analyze the parameters $n_1, n_2$ and $n_3$ associated with $D$ such that $n_1 (sk)= n_2 k + n_3 (k+1)$. 
Clearly, $n_2=\sum_{j\in\Z_k} n_2^j$, $n_3=k$ and then $n_1s=n_2+k+1$. 

In order to verify  item 2. of Theorem \ref{minorsgrafo}, it only remains to prove that $n_1\leq k-1$ which is equivalent to prove that $n_2\leq k(s-1)-2$.

From the definition, $n_2^j\leq s-1$ for every $j\in\Z_k$.

Suppose that $n_2^j=s-1$ for every $j\in \Z_k$, $j\neq \ell$.

If $\ell=k-1$, then, $t_j=t_{j+1}$ for all $0\leq j\leq k-2$ and then $t_0=t_j$ for all $1\leq j\leq k-1$. In this case, $C^{i_0}=W$ contradicting the fact that $|C^{i_0}\cap \overline{W}|\geq 1$.

Now, if $\ell \neq k-1$ then it holds that  $t_j=t_{\ell}$ for all $0\leq j\leq \ell$ and $t_{j}=t_{\ell +1}$ for all $\ell+1\leq j\leq k-1$. Since $n_2^{k-1}=s-1$ then either $t_{k-1}=t_0-1$ and $t_0\neq 0$ or $t_{k-1}=s-1$ and $t_0=0$. In any case $C^{i_{\ell}}=W$ again contradicting $|C^{i_{\ell}}\cap \overline{W}|\geq 1$.

As a consequence, there must be at least two different values of $j$ with $n_2^j\leq s-2$ and the theorem follows.
\qed

\begin{corollary}
Every facet defining $(s+1)$-inequality of $Q(C_{sk}^k)$ is a relevant minor inequality verifying Conjecture \ref{conj2}. 
\end{corollary}

\proof

From the proof of the previous theorem every facet defining $(s+1)$-inequality of $Q(C^k_{sk})$ is associated with a relevant minor $C^{k'}_{n'}$ with $n'=sk-(n_2+k)$, $k'=k-n_1$ and $n_1s=n_2+k+1$.
 
Then, $n'=sk-(n_1s-1)=s(k-n_1)+1=sk'+1$ and Conjecture \ref{conj2} holds.
\qed

\bigskip

In the remainder of this section we will prove that the $(s+1)$-inequalities in the form (\ref{minor s+1})
with $W$ such that $|x^i\cap W|=1$ for every $i\in \Z_k$ can be separated in polynomial time.

\medskip

Any of these inequalities can be written as 
$$\sum_{i\in W} x_i+\sum_{i\in\Z_{sk}} x_i  \geq  s+1$$ 
or equivalently
\begin{equation}\label{facetamin}
\sum_{i\in W} x_i  \geq  s+1- \sum_{i\in\Z_{sk}} x_i.
\end{equation}

Defining $L(x):= s+1- \sum_{i\in\Z_{sk}} x_i$, the separation problem for these inequalities can be stated as follows: given $\hat{x} \in \mathbb{R}^n$, decide if there exists $W\subset \Z_{sk}$ 
with $|x^i\cap W|=1$ for all $i\in \Z_k$ such that 
$$\sum_{i\in W} \hat{x}_i  <  L(\hat{x}).$$

We will reduce this problem to a shortest path problem in an acyclic digraph.

\begin{figure}[h]
	\centering
		\includegraphics[width=0.9\textwidth]{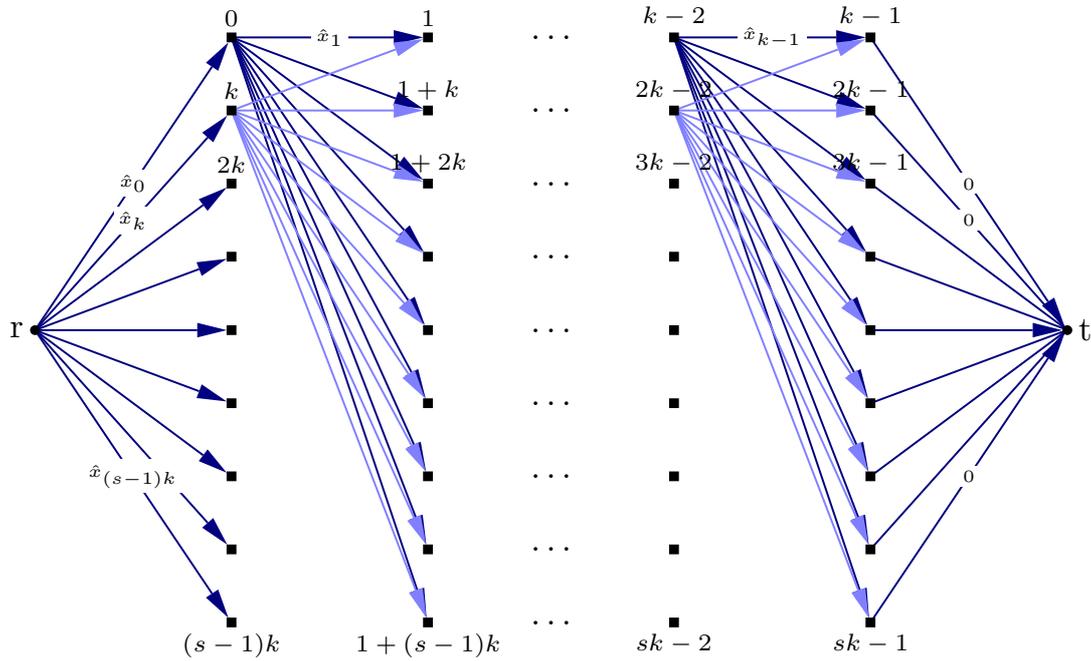}
	\caption{The digraph $D(C_{sk}^k)$.}
	\label{fig:digrafo}
\end{figure}

For this purpose let us define the digraph $D(C_{sk}^k)=(V,A)$, where 
$$V=\left(\bigcup_{i\in\Z_k} V_i\right) \cup \{r,t\}$$ with $V_i=x^i=\{i,i+k,\ldots,i+(s-1)k\}$ for $i\in\Z_k$ and 

$$A=\left(\bigcup_{i\in\Z_{k-1}} A_i\right)\cup A_r \cup A_t$$
with $A_i=\{(l,m):l\in V_i,\;  m\in V_{i+1}\}$ for $i \in \Z_{k-1}$, $A_r=\{(r,m): m\in V_0\}$ and $A_t=\{(l,t):l\in V_{k-1}\}$. For illustration see Figure \ref{fig:digrafo}.

Observe that there is a one-to-one correspondence between $r t$-paths in $D(C_{sk}^k)$ and subsets $W\subset \Z_{sk}$ with $|x^i\cap W|=1$ for every $i\in \Z_{k}$.
 
Indeed, let $W\subset \Z_{sk}$ such that $x^j\cap W=\{i_j\}$ for every $j\in \Z_{k}$. Clearly, $\{r,i_0,i_1,\dots,i_{k-1},t\}$ induces an $r t$-path in $D(C_{sk}^k)$.
Conversely, if $P$ is an $r t$-path in $D(C_{sk}^k)$, by construction, 
$|V(P)\cap V_i|=1$ for all $i\in\Z_k$, and 
$$W=\bigcup_{i\in\Z_k} (V(P)\cap V_i)$$ 
verifies $|W\cap x^i|=1$ for every $i\in\Z_k$.

Then, we have the following result.

\begin{theorem}\label{poly}
Given $s\geq 2$ and $k \geq 4$, the inequalities (\ref{minor s+1}) with $|x^i\cap W|=1$ for every $i\in \Z_{k}$, can be separated in polynomial time.
\end{theorem}

\proof 
Let $\hat{x} \in \mathbb{R}^n$ and $D(C^k_{sk})$ be the digraph previously defined. For each arc $(i,j)\in A\setminus A_t$, assign length $\hat{x}_j$ otherwise if $(i,j)\in A_t$ assign length zero (see Figure \ref{fig:digrafo}). 

In this way, given $W\subset \Z_{sk}$ with $|W\cap  x^i|=1$,  
the length of the $r t$-path corresponding to $W$ is $\sum_{i\in W} \hat{x}_i$.

Then, the separation problem can be reduced to decide if there exists an $rt$-path in  $D(C_{sk}^k)$ with length less than $L(\hat{x})=s+1- \sum_{i\in \Z_{sk}} \hat x_i$ or, equivalently, if the shortest path in $D(C_{sk}^k)$ has length less than $L(\hat{x})$. 

Since $D(C_{sk}^k)$ is acyclic, computing the shortest path can be done in polynomial time using for instance Bellman algorithm \cite{bel}.
\qed

\medskip

Let us observe that, if $\hat{x}\in Q_R(C_{sk}^k)$, the separation problem for inequalities of the form (\ref{minor s+1}) such that $|x^i\cap W|=1$ is equivalent to the separation problem for relevant minor $(s+1)$-inequalities of $Q(C_{sk}^k)$.
%
\section{The Set Covering Polyhedron of $C^k_{2k}$ and  $C^k_{3k}$}

In this section we prove that minor inequalities together with the boolean facets completely describe the set covering polyhedron of $C^k_{2k}$ and  $C^k_{3k}$ for $k\geq 2$. Actually, if $k=2$ the matrices are ideal and if $k=3$ the result follows from \cite{Mah} and Theorem \ref{minor}. 

The key is to prove that for every non boolean non rank facet defining inequality 
there exists a cover $\tilde{x}$ of $C^{k}_{sk}$, with $|\tilde{x}|=s+1$ and $\tilde{x}\cap W=\emptyset$. 
Provided such a cover exists, then,   
$$\alpha \leq \sum_{i\in W} a_i \tilde{x}_i + \sum_{i\in \overline{W}} a^0 \tilde{x}_i= a^0 |\tilde{x}|= a^0 (s+1).$$ 
Hence, $\alpha= (s+1)a^0$ and  (\ref{desig}) is a facet defining $(s+1)$-inequality. By Theorem \ref{minor}, it is a minor inequality. 

Recall that if (\ref{desig}) is a facet defining inequality of $Q(C^{k}_{sk})$ 
for some $W\subset \Z_{sk}$ and $a^0\geq 1$, $|x^i \cap W|\geq 1$ for all $i\in \Z_k$ and, from Theorem \ref{wraya} item 1, $|C^{i}\cap \overline{W}|\geq 2$, for all $i\in \Z_k$.

Let us start with the case $s=2$.

\begin{theorem}\label{vale2}
For every $k\geq 3$, every non boolean facet defining inequality of $Q(C^k_{2k})$ is a minor inequality.
\end{theorem}

\proof
Let $W\subset \Z_{2k}$ such that (\ref{desig}) is a facet defining inequality of $Q(C^{k}_{2k})$.
As we have already observed, it is enough to prove that there is a cover $\tilde{x}$ such that $|\tilde{x}|=3$ and $\tilde{x}\subset \overline{W}$.  
W.l.o.g we can assume that $0\in \overline W$ and since $|x^0\cap W|\geq 1$ we have that $k\in W$.

Considering that $|C^{0}\cap \overline{W}|\geq 2$, we can define  $t=\max \,\{s: 1\leq s \leq k-1\,,\; s\in \overline{W}\}$
that makes 
$(t,k]_{2k}\subset W$. Since $|W\cap x^t|\geq 1$ and $t\in \overline W$, $t+k\in W$.

Since $|C^t \cap \overline W|\geq 2$ and $(t,k]_{2k}\subset W$, there exists $t'\in  \overline W\cap (k,t+k)_{2k}$.
We have that $\tilde{x}=\{0, t, t'\}\subset \overline W$ is a cover of $C^k_{2k}$ and the proof is complete. 
\qed

\medskip

In what follows we prove that for every inequality (\ref{desig}) defining a facet of $Q(C^k_{3k})$ there exists a cover $\tilde{x}$ of $C^k_{3k}$ with $|\tilde{x}|=4$ and $\tilde{x}\cap W=\emptyset$. For this purpose, let us call $\mathcal{W}$ the family of $W\subset \Z_{3k}$ such that, for every $i\in \Z_{3k}$, $|x^i\cap W|\geq 1$ and $|C^i \cap \overline W|\geq 2$. Clearly, every $W$ associated with a non boolean facet defining inequality of $Q(C^k_{3k})$ is in $\mathcal{W}$.   

For every $W\in \mathcal{W}$ and $i\in \overline{W}$ we define 
$$\omega(i):=\min\{t\geq 0: i+k+t\in\overline{W}\}.$$ 
Since $|C^{i+k}\cap \overline W|\geq 2$, $\omega(i)\leq k-2$ for all $i\in \overline{W}$.

Firstly, we have:

\begin{lemma}\label{lemacaso1}
Let $W\in \mathcal{W}$ and  $i\in \overline W$ such that 
$\overline{W}\cap [i+2k,i+2k+\omega(i)]_{3k}\neq \emptyset$. Then there exists a cover $\tilde{x}$ of $C^k_{3k}$ with $|\tilde{x}|=4$ and $\tilde{x}\cap W=\emptyset$.
\end{lemma}

\proof
Let us first observe that $\omega(i)\geq 1$. Indeed, if $\omega(i)=0$, $i+k\in \overline{W}$ and $[i+2k,i+2k+\omega(i)]_n =\{i+2k\}$. Then, by hypothesis, $i+2k \in \overline W$ contradicting the fact that $|x^i\cap W|\geq 1$. 

Let $\ell \in \overline{W}\cap [i+2k,i+2k+\omega(i)]_{3k}$. Since $|C^{i+w(i)}\cap \overline{W}|\geq 2$ and  $[i+k,i+k+\omega(i))_{3k}\subset W$, there exists $t\in \overline{W} \cap[i+w(i),i+k)_{3k}$.
Then, 
$\tilde{x}=\{i, t, i+k+\omega(i), \ell \}\subset \overline{W} $ is a cover of $C^k_{3k}$.
\qed

\medskip

As immediate consequence we get that if $W\in \mathcal{W}$ and $|W\cap x^j|=1$ for some $j\in \Z_{3k}$ then $C^k_{3k}$ admits a cover $\tilde{x}\subset \overline W$ with $|\tilde{x}|=4$. Indeed, w.l.o.g. we can assume that $j\in \overline W$ and $j+k\in W$. Then, it holds that  
$\overline{W}\cap [j+2k,j+2k+\omega(j)]_{3k}\neq \emptyset$ and we can apply the above lemma. 

Therefore, it only remains to consider subsets $W\in \mathcal W$ such that, for every $i\in\Z_{3k}$, $|x^i \cap W|\geq 2$  and for every $i\in \overline{W}$, $[i+2k,i+2k+\omega(i)]_{3k}\subset W$. W.l.o.g. we can assume that $0\in\overline{W}$. Let us call $\mathcal{W}^*$ the family of all subsets $W\in \mathcal W$ satisfying these conditions. 

Given $W\in \mathcal{W}$ and $i\in \overline{W}$ let us recursively define the sequence $r^i= \{r^i_t\}_{t=0}^{\infty}\subset \Z_{3k}$  as follows:
\begin{itemize}
\item $r^i_0=i$
\item $r^i_t= r^i_{t-1}+k+ \omega(r^i_{t-1})$, for $t\geq 1$. 
\end{itemize}
According to the sequence $r^i$, we define $p^i=\max\left\{t: \sum_{j=0}^{t-1} \omega(r^i_j) \leq k-1\right\}$.

Observe that, given $i\in \overline{W}$, by definition of $\omega(i)$,  $r^i\subset \overline{W}$ and 
$[r^i_t+k, r^i_{t+1})_{3k}\subset W$ for all $t$. Moreover, if $W\in \mathcal{W}^*$, $[r^i_t+2k, r^i_t+2k+\omega(r^i_t)]_{3k}=[r^i_t+2k, r^i_{t+1}+k]_{3k}\subset W$ and, since $[r^i_{t+1}+k, r^i_{t+2})_{3k}\subset W$, we have that $[r^i_t+2k, r^i_{t+2})_{3k}\subset W$. 

From this observation we can conclude that for every $W\in \mathcal{W}^*$ and every $i\in \overline{W}$, $p^i\neq 0$ (mod $k$).

\begin{figure}[h]
	\centering
		\includegraphics[width=0.75\textwidth]{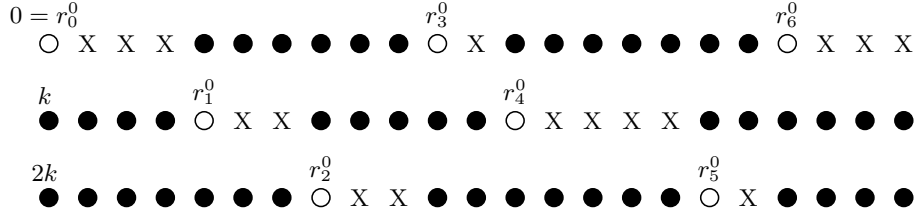}
			\caption{A subset $\widetilde{W}\in \mathcal {W}$.}
	\label{fig:escaleras}
\end{figure}

In Figure \ref{fig:escaleras} we sketch a subset $\widetilde{W}\in \mathcal W$, for $k=23$.
Each node corresponds to an element in $\Z_{69}$. Black nodes and white nodes correspond to elements in and out $\widetilde{W}$, respectively. Crosses correspond to elements that may or may not belong to $\widetilde{W}$. We also show the first seven elements of the sequence $r^0$ that allow us to see that $p^0=6$, i.e. $p^0=0$ (mod 3). 
We will show that $\widetilde{W}\notin \mathcal{W}^*$.

In fact, observe that $r^0_7=r^0_2$ and the sequence $r^0$ cycles. Then, $r^0_0\in [r^0_6+2k, r^0_7+k)_{69}$ and $r^0_0\notin \widetilde{W}$. Since for all $W\in \mathcal{W}^*$ and every $t$, $[r^i_t+2k, r^i_{t+1}+k]_{3k}\subset W$, we have that  $\widetilde{W}\notin \mathcal{W}^*$.

In general, given $W\in \mathcal W$ such that $p^i=0$ (mod 3) for some $i\in \overline W$, following the same reasoning we can arrive to $r^i_{p^i+1}=r^i_2$ and $r^i_0\in \overline W \cap [r^i_{p^i}+2k, r^i_{p^i+1}+k)_{3k}$ 
and then $W\notin \mathcal {W}^*$. Then, we have:

\begin{lemma}
For every $W\in \mathcal{W}^*$, $p^i\neq 0$ (mod 3) for every $i\in \overline{W}$
\end{lemma}

Moreover, we have the following result:

\begin{lemma}
Let $W\in \mathcal{W}^*$. Then, there exists $i\in \overline W$ such that $p^i=1$ \emph{(mod 3)}.
\end{lemma}

\proof
Let $i\in \overline W$ and let us call $j=r^i_{p^i}$. 

\begin{figure}[h]
	\centering
		\includegraphics[width=0.65\textwidth]{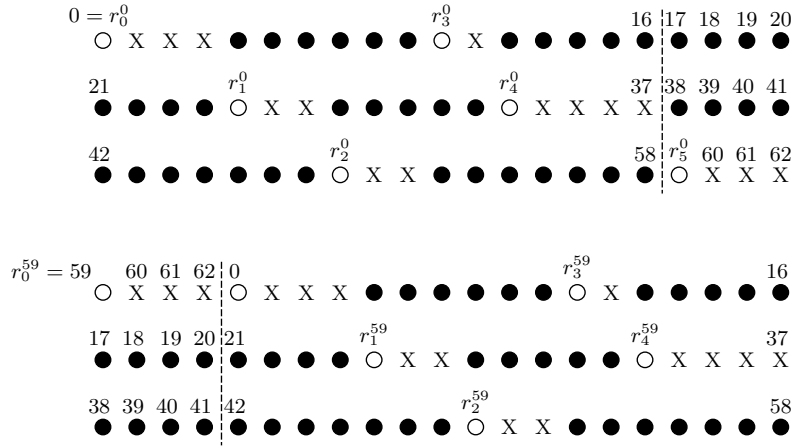}
	\caption{$W\in \mathcal{W^*}$ with $p^0=5=2$ (mod 3), $j=r^0_{p^0}=59$ and $p^{59}=4=1$ (mod 3).}
	\label{rW5}
\end{figure}

According to the previous lemma, 
$p^i\neq 0$ (mod 3) and then, consider 
$p^i= 2$ (mod 3). W.l.o.g. we can assume that $i=0$ (for illustration consider the example in Figure \ref{rW5}). 

We have that $[j+k,k)_{3k}\subset W$. In addition, we know that $[r^0+k,r^0_1)_{3k}=[k,r^0_1)_{3k}\subset W$.
Then, $[j+k,r^0_1)_{3k}\subset W$ implying that $r^j_1=r^0_1$ and $r^j_t=r^0_t$ for all $t\geq 1$. Moreover,  
$$\omega(r^j_0)=\omega(j)= k- \sum_{t=0}^{p^0-1} \omega(r^{0}_t)+ \omega(r^0_0).$$

Then, 
$$\sum_{t=0}^{p^0-2} \omega(r^{j}_t)\leq k-1 \;\; \mathrm{ and } \;\; \sum_{t=0}^{p^0-1} \omega(r^{j}_t)> k-1.$$
Therefore, $p^j=p^0-1$ and $p^j=1$ (mod 3).
\qed

Finally, we can prove:

\begin{theorem}\label{lemacaso2}
Let $W\in \mathcal{W^*}$ associated with a facet defining inequality (\ref{desig}). Then, there exists a cover $\tilde{x}$ of $C^k_{3k}$ with $|\tilde{x}|=4$ and $\tilde{x}\cap W=\emptyset$.
\end{theorem}

\proof
By the previous observation, w.l.o.g. we can assume that $p^0=1$ (mod 3). 
It is not hard to see that if $p^0=1$, $C^{2k}\subset W$, a contradiction. Then, $p^0\geq 4$. We will prove that  $p^0\neq 4$.
 
\begin{figure}[h]
	\centering
		\includegraphics[width=0.65\textwidth]{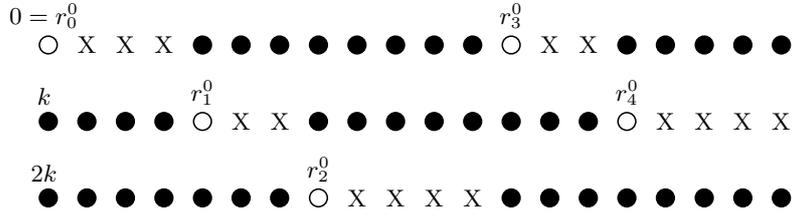}
		\caption{A subset $W\in \mathcal {W}^*$ with $p^0=4$.}
	\label{fig:rW4}
\end{figure}

Suppose that $p^0= 4$ (see Figure \ref{fig:rW4} as example).

In this case, the set  $\{r_0^0,r^0_3,r^0_1,r^0_4,r^0_2\}\subset \overline{W}$ results a cover of $C^k_{3k}$ with cardinality five, $\alpha \leq 5a^0$ and every root of (\ref{desig}) has cardinality four or five.

From Lemma \ref{prop} item 1.(c), there exists a root $\tilde{x}$ of (\ref{desig}) such that $\tilde{x} \cap C^{2k}=\{s, t\}$ with $s\in [2k,t)_{3k}$. We will prove that $\{s,t\}\subset \overline{W}$.

Observe that $C^{2k}\cap \overline{W}\subset [r^0_2,r^0_3+2k)_{3k}$.
If $t \in W$, applying Lemma \ref{prop} item 2.(a), we have that $0\in W$, a contradiction. 
Hence, $t \in \overline{W}$.

Again applying Lemma \ref{prop} item 2.(a), if $s \in W$ then $[t-k,2k)_{3k}\cap C^{k}\subset W$ implying $r^0_4\in W$, a contradiction. Then, $s \in \overline{W}$. 

Then, $\{s,t\}\subset [r^0_2,r^0_2+\omega(r^0_2))_{3k}$.

Since $\tilde{x}$ is a minimal cover, there exist 
$\ell, \ell' \in \tilde{x}$ such that $\ell\in [s-k,t-k)\subset W$ and $\ell' \in (s+k,t+k]\subset W$. 

Recall that $|\tilde{x}|=4$ or $|\tilde{x}|=5$. Moreover, $|\tilde{x}|=5$ if and only if $\tilde{x}\subset \overline W$.
Then, $|\tilde{x}|=4$ and $\tilde{x}=\{s,t,\ell,\ell'\}$.

Since $\tilde{x}$ is a root of (\ref{desig}), $\alpha = 2a^0+a_{\ell}+a_{\ell'}$.

Observe that $\hat{x}=\{r^0_1,r^0_4,\ell',t\}$ is a cover of $C^k_{3k}$ which violates (\ref{desig}) since 
$3a^0+a_{\ell'} < \alpha = 2a^0+a_{\ell}+a_{\ell'}$. 

Therefore, $p^0\neq 4$ and then $p^0\geq 7$.


Finally, 
$\tilde{x}=\{r^0_0, r^0_6,r^0_4, r^0_2\}\subset \overline W$ is a cover of $C^k_{3k}$.
\qed

\medskip

Summarizing, we have proved the following: 

\begin{theorem} \label{vale}
For every $k\geq 3$, every non boolean facet defining inequality of $Q(C^k_{3k})$ is a minor inequality.
\end{theorem}

\proof
Let $W\subset \Z_{3k}$ such that (\ref{desig}) is a facet defining inequality of $Q(C^{k}_{3k})$ for some $a^0\geq 1$.

As we have already observed, it is enough to prove that there exists a cover $\tilde{x}$ such that $|\tilde{x}|=4$ and $\tilde{x}\subset \overline{W}$.  

Therefore, by Lemmas \ref{lemacaso1} and \ref{lemacaso2}, the theorem follows.
\qed

\bigskip

Let us observe that, from Theorems \ref{poly}, \ref{vale2} and \ref{vale} we have a polynomial time algorithm based on linear programming that solves the Set Covering Problem for matrices $C^k_{2k}$ and $C^k_{3k}$, for every $k\geq 3$.

On the other hand, Proposition 5.5 in \cite{Nes} proves that, for every $k$, the set covering polyhedron of $C^k_{2k+1}$, $C^k_{3k+1}$, $C^k_{3k+2}$ is described by means of boolean facets and the rank constraint. Lemma \ref{menor_circ} and Theorems \ref{vale2} and  \ref{vale} give an alternative proof of the same fact. 

Indeed, remind that the rank constraint is always a facet defining inequality of $Q(C_{sk+r}^k)$ when $1\leq r\leq k-1$ and it has right hand side $s+1$. Hence, if $Q(C_{sk+r}^k)$ with $s=2,3$ and $1\leq r\leq s-1$ has a non boolean non rank facet defining inequality of the form (\ref{desig}) it must have right hand side at least $s+2$. But, applying Lemma \ref{menor_circ}, $C_{sk+r}^k$ is a minor of $C^{k'}_{sk'}$ for some $k'$ and then there would exist some facet defining inequality for $Q(C^{k'}_{sk'})$ with right hand side different from $s+1$ contradicting Theorems \ref{vale2} and \ref{vale}.

\end{document}